\documentclass[leqno]{article}

\usepackage{amssymb,amsmath,latexsym}

\newtheorem{thm}{Theorem}[section]
\newtheorem{pro}[thm]{Proposition}
\newtheorem{cor}[thm]{Corollary}
\newtheorem{lem}[thm]{Lemma}

\newcommand{\loc}{\mathrm{loc}}
\newcommand{\Ker}{\mathrm{Ker}}
\newcommand{\im}{\mathrm{Im}}
\newcommand{\Hom}{\mathrm{Hom}}
\newcommand{\End}{\mathrm{End}}
\newcommand{\ext}{\mathrm{Ext}}
\newcommand{\ann}{\mathrm{Ann}}
\newcommand{\irr}{\mathrm{Irr}}
\newcommand{\soc}{\mathrm{Soc}}
\newcommand{\mat}{\mathrm{Mat}}
\newcommand{\sol}{\mathrm{Solv}}
\newcommand{\gl}{\mathrm{gl.dim}}
\newcommand{\id}{\mathrm{inj.dim}}
\newcommand{\hoe}{\mathrm{ht}}
\newcommand{\ug}{\mathrm{u.gr}}
\newcommand{\lgr}{\mathrm{l.gr}}

\newcommand{\ag}{\mathfrak{a}}
\newcommand{\s}{\mathfrak{s}}
\newcommand{\slg}{\mathfrak{sl}}

\newcommand{\nb}{\mathbb{N}}
\newcommand{\cb}{\mathbb{C}}
\newcommand{\fb}{\mathbb{F}}
\newcommand{\db}{\mathbb{D}}

\newcommand{\oc}{\mathcal{O}}
\newcommand{\ic}{\mathcal{I}}
\newcommand{\pc}{\mathcal{P}}

\newcommand{\mc}{\mathcal{M}}
\newcommand{\nc}{\mathcal{N}}
\newcommand{\ac}{\mathcal{A}}

\begin{document}


\title{\bf Injective Modules and Prime Ideals of Universal Enveloping Algebras}

\author{J\"org Feldvoss\thanks{E-mail address: \tt jfeldvoss@jaguar1.usouthal.edu}\\
{\small Department of Mathematics and Statistics}\\{\small University of South Alabama}\\
{\small Mobile, AL 36688--0002, USA}}

\date{Dedicated to Professor Edgar E. Enochs on the occasion of his 72$^{\rm nd}$ birthday}

\maketitle


\begin{abstract}

\noindent In this paper we study injective modules over universal enveloping algebras 
of finite-dimensional Lie algebras over fields of arbitrary characteristic. Most of
our results are dealing with fields of prime characteristic but we also elaborate on 
some of their analogues for solvable Lie algebras over fields of characteristic zero.
It turns out that analogous results in both cases are often quite similar and resemble 
those familiar from commutative ring theory.   
\medskip

\end{abstract}


\section*{Introduction}

\noindent In this paper we investigate the injective modules and their relation to 
prime ideals in universal enveloping algebras of finite-dimensional Lie algebras. 
Especially, in the case that the ground field is of prime characteristic we obtain 
several results that seem to be new. It should be remarked that most of the results 
of the first two sections and the last section are already contained in an unpublished 
manuscript of the author (cf.~\cite{Fel90}) but the entire third section and Theorem
\ref{lastermininjres} are completely new. In the following we will describe the contents 
of the paper in more detail.

The first section provides the framework for the paper. We begin by recalling the
well-known result from noetherian ring theory that every injective module decomposes 
uniquely (up to isomorphism and order of occurrence) into a direct sum of indecomposable 
injective modules. Then it is shown that universal enveloping algebras of finite-dimensional 
Lie algebras over fields of prime characteristic are FBN rings. As a consequence, 
indecomposable injective modules are in bijection with prime ideals. Moreover, it is 
proved that the universal enveloping algebra of a finite-dimensional Lie algebra 
over a field of prime characteristic is a Matlis ring (i.e., every indecomposable 
injective module is the injective hull of a prime factor ring of the universal 
enveloping algebra considered as a one-sided module) if and only if the underlying 
Lie algebra is abelian. A similar result might also hold in characteristic zero but 
we were neither able to prove this nor to find it in the literature.

In the second section we study certain finiteness conditions for injective hulls. 
It is well-known from a result obtained by Donkin \cite{Don82} and independently by
Dahlberg \cite{Dah84} that injective hulls of locally finite modules over universal 
enveloping algebras of finite-dimensional solvable Lie algebras over fields of 
characteristic zero are again locally finite. We show that the converse of this result 
holds, i.e., the locally finiteness of injective hulls of locally finite modules in 
characteristic zero implies that the underlying Lie algebra is solvable. In fact, the 
locally finiteness of the injective hull of the one-dimensional trivial module already 
implies that the underlying Lie algebra is solvable. This generalizes an observation of 
Donkin in \cite{Don82}. Moreover, we prove that every essential extension of a locally 
finite module over a universal enveloping algebra of any finite-dimensional Lie algebra 
over a field of prime characteristic is locally finite by applying a result of Jategaonkar 
\cite{Jat74} in conjunction with the result from the first section saying that universal 
enveloping algebras of finite-dimensional Lie algebras over fields of prime characteristic 
are FBN rings. In particular, injective hulls of locally finite modules are always locally 
finite. By generalizing slightly another result of Jategaonkar \cite{Jat75}, we also 
show that injective hulls of artinian modules over universal enveloping algebras of 
finite-dimensional Lie algebras over a field of prime characteristic are always artinian. 
Finally, it is established that for the universal enveloping algebra of a non-zero 
finite-dimensional Lie algebra over a field of prime characteristic non-zero noetherian 
modules are never injective by proving that the injective dimension of a non-zero noetherian 
module coincides with the dimension of the underlying Lie algebra. On the other hand, there 
are artinian and locally finite modules of any possible injective dimension.   

In the third section we consider certain locally finite submodules of the linear dual of 
a universal enveloping algebra. We start off by showing how an argument from \cite{Dah84} 
can be changed slightly to make it work over arbitrary fields of any characteristic and 
therefore obtaining a different (and in our opinion more transparent) proof of a result due 
to Levasseur \cite{Lev76}. Then we give a very short proof of the main result of \cite{Kos54}
by using the locally finiteness of injective hulls of locally finite modules over universal 
enveloping algebras of finite-dimensional solvable Lie algebras in characteristic zero in 
an essential way. In fact, this argument was motivated by our proof of the injectivity of 
the continuous dual of the universal enveloping algebra of an arbitrary finite-dimensional 
Lie algebra over a field of prime characteristic. As an immediate consequence, we obtain 
that in prime characteristic the cohomology with values in the continuous dual vanishes 
in every positive degree. In particular, Koszul's cohomological vanishing theorem does 
remain valid in prime characteristic. These results seem to be new. Moreover, the modular 
cohomological vanishing theorem is much stronger than its analogue in characteristic zero 
which follows from a recent result of Schneider (cf.~\cite{Mas00}) and says that the 
cohomology with values in the continuous dual vanishes in degrees one and two.

The last section closes the circle of ideas by coming back to the correspondence
between injective modules and prime ideals. It is verified that universal enveloping 
algebras of finite-dimensional Lie algebras over fields of prime characteristic 
are injectively homogeneous in the sense of \cite{BH88}. As a consequence of the 
general theory of injectively homogeneous rings developed in \cite{BH88} we obtain 
a nice description of a minimal injective resolution of the universal enveloping 
algebra as a module over itself in terms of the injective hulls of its prime factor 
rings considered as one-sided modules. In particular, this enables us to show that 
the last term of such a minimal injective resolution is isomorphic to the continuous 
dual which was proved by Barou and Malliavin \cite{BM85} for finite-dimensional solvable 
Lie algebras over algebraically fields of characteristic zero.
 
Throughout this paper we will assume that all associative rings have a unity element 
and that all modules over associative rings are unital. 


\section{Injective Modules and Prime Ideals}


\noindent Since the universal enveloping algebra of a finite-dimensional Lie algebra 
$\ag$ is left and right noetherian (cf.~\cite[Theorem V.6]{Jac79}), finding all injective 
left and right $U(\ag)$-modules reduces to the classification of the indecomposable 
ones (see \cite[Theorem 2.5, Proposition 2.6, and Proposition 2.7]{Mat58}):

\begin{pro}
Let $\ag$ be a finite-dimensional Lie algebra over an arbitrary field. Then the 
following statements hold:
\begin{enumerate}
\item[{\rm (1)}] Every injective left or right $U(\ag)$-module is a direct sum of 
                 indecomposable injective submodules.
\item[{\rm (2)}] If I is an indecomposable injective left or right $U(\ag)$-module, 
                 then $\End_\ag(I)$ is local. In particular, the decomposition in the 
                 first part is unique up to isomorphism and order of occurrence of the 
                 direct summands.\quad $\Box$
\end{enumerate}
\end{pro}

\noindent In order to be able to parameterize the indecomposable injective left or right
$U(\ag)$-modules, one needs the following concept from non-commutative ring theory. A 
left and right noetherian associative ring $R$ is called a {\it FBN ring\/} if every 
essential left ideal and every essential right ideal of every prime factor ring of $R$ 
contains a non-zero two-sided ideal (which, in fact, is essential). While classifying 
the indecomposable injective $U(\ag)$-modules by analogy with the commutative case 
(see \cite[Proposition 3.1]{Mat58}), one should be aware that the injective hull of 
$U(\ag)/\pc$ (considered as a left or right $U(\ag)$-module) is not necessarily 
indecomposable for every prime ideal $\pc$ of $U(\ag)$. For example, the injective 
hull of $U(\ag)/\ann_{U(\ag)}(S)$ is isomorphic to the direct sum of $d$ copies of 
the injective hull of any simple $\ag$-module $S$ of dimension $d>1$ (cf.~the proof 
of Theorem \ref{mat} and Theorem \ref{lastermininjres}).

Let $M$ be a non-zero $U(\ag)$-module. A two-sided ideal $\pc$ is said to be {\it 
associated to\/} $M$ if there exists a submodule $N$ of $M$ such that $\pc$ equals 
the annihilator of every non-zero submodule of $N$. It is well-known that $\pc$ is 
necessarily prime and that for an indecomposable injective module $I$ there exists 
a unique prime ideal $\pc_I$ associated to $I$ (cf.~\cite{BGR73}).

If $\ag$ is a {\it finite-dimensional\/} Lie algebra over a field of {\it prime\/}
characteristic, then $U(\ag)$ is a {\it finitely generated\/} $C(U(\ag))$-module 
(cf.~\cite[Theorem 5.1.2]{SF88}). Hence one has the following well-known facts which 
are crucial for the results obtained in this paper:

\begin{itemize}
\item[(IC)] $U(\ag)$ is {\it integral\/} over its center $C(U(\ag))$ (cf.~\cite[Theorem 
            6.1.4]{SF88}). More generally, there exists a subalgebra $\oc(\ag)\cong\fb
            [t_1,\dots,t_{\dim_\fb\ag}]$ of $C(U(\ag))$ such that $U(\ag)$ is {\it 
            integral\/} over every subring $C$ of $U(\ag)$ with $\oc(\ag)\subseteq C
            \subseteq C(U(\ag))$.
\item[(PI)] $U(\ag)$ is a {\it PI ring\/} (cf.~\cite[p.~xi]{GW89}).
\end{itemize}
\smallskip

The next result shows that the indecomposable injective modules over universal
enveloping algebras in prime characteristic can be classified by their associated
prime ideals.

\begin{thm}\label{FBN}
Let $\ag$ be a finite-dimensional Lie algebra over a field of prime characteristic. 
Then the universal enveloping algebra $U(\ag)$ is a FBN ring. In particular, there 
is a one-to-one correspondence between the indecomposable injective $U(\ag)$-modules 
and the prime ideals of $U(\ag)$ given by $I\mapsto\pc_I$, where $\pc_I$ is the unique 
prime ideal associated to $I$.
\end{thm}

\noindent {\it Proof.\/} The first assertion follows from \cite[Proposition 8.1(b)]{GW89} 
and the second assertion is a consequence of the first and \cite[Theorem 3.5]{Kra72}.
\quad $\Box$
\bigskip

\noindent {\bf Question.} Let $\ag$ be a finite-dimensional Lie algebra over a 
field of characteristic {\it zero\/}. It would be interesting to know when $U(\ag)$ 
is a FBN ring? Is $U(\ag)$ only an FBN ring if $\ag$ is {\it abelian\/}?
\bigskip

An associative ring $R$ is called a {\it Matlis ring\/} if every indecomposable injective 
left or right $R$-module is isomorphic to the injective hull of $R/\pc$ (considered as a 
left or right $R$-module) for some prime ideal $\pc$ of $R$. Every left and right noetherian 
Matlis ring is a FBN ring (see \cite[Corollary 3.6]{Kra72}), but the converse is not true as 
follows from Theorem \ref{FBN} and the next result. 

\begin{thm}\label{mat}
Let $\ag$ be a finite-dimensional Lie algebra over a field of prime characteristic. Then 
the universal enveloping algebra $U(\ag)$ is a Matlis ring if and only if $\ag$ is abelian.
\end{thm}

\noindent {\it Proof.\/} Since both conditions are independent of the ground field $\fb$, 
we can assume that $\fb$ is algebraically closed. Suppose that $U(\ag)$ is a Matlis ring. 
According to \cite[Corollary 14]{Kra70}, every prime ideal of $U(\ag)$ is completely prime. 
Let $S$ be a simple $\ag$-module and set $\db:=\End_\ag(S)$. Since $S$ is finite-dimensional
(cf.~\cite[Theorem 5.2.4]{SF88}), $\db$ is a finite-dimensional division algebra over $\fb$, 
and thus $\db=\fb$. Then the density theorem (cf.~\cite[Theorem 16, p.~95]{Kap72}) implies 
that $$U(\ag)/\ann_{U(\ag)}(S)\cong\End_\fb(S)\cong\mat_d(\fb),$$ where $d:=\dim_\fb S$. 
Since $S$ is simple, $\ann_{U(\ag)}(S)$ is primitive (i.e., prime), and thus, $\ann_{U(\ag)}
(S)$ is completely prime. It follows that $\mat_d(\fb)$ has no zero divisors, i.e., $d=1$. 
Hence every simple $\ag$-module is one-dimensional. By virtue of a result due to Jacobson, 
there exists a (finite-dimensional) faithful semisimple $\ag$-module (see \cite[Theorem 
5.5.2]{SF88}). Therefore, we have $$[\ag,\ag]\subseteq\bigcap_{S\in\irr(\ag)}\ann_\ag(S)=0,$$ 
where $\irr(\ag)$ denotes the set of isomorphism classes of simple $\ag$-modules, i.e., $\ag$ 
is abelian. Finally, the converse is just \cite[Proposition 3.1]{Mat58}.\quad $\Box$
\bigskip

\noindent {\bf Remark.} The proof of Theorem \ref{mat} applied to a composition factor $S$ 
of the adjoint module of a finite-dimensional Lie algebra $\ag$ over a field of characteristic 
zero shows that in this case the universal enveloping algebra $U(\ag)$ can only be a Matlis 
ring if $\ag$ is {\it solvable\/} (cf.~also \cite[p.~49]{BGR73}). This still leaves the question
as to whether Theorem \ref{mat} is also true in characteristic zero.


\section{Injective Hulls}


\noindent In this section several finiteness properties of injective hulls are considered. Let 
$R$ be an associative ring and let $M$ be a left or right $R$-module. An injective module $I$ 
is called an {\it injective hull\/} (or an {\it injective envelope\/}) of $M$ if there exists 
an $R$-module monomorphism $\iota:M\to I$ such that the image $\im(\iota)$ of $\iota$ is an 
essential submodule of $I$. (By abuse of language, the pair $(I,\iota)$ is also called an
injective hull of $M$.)

It is well-known that every module has an injective hull (cf.~\cite[Theorem 4.8(a)]{GW89}). 
Moreover, injective hulls satisfy the following {\it universal properties\/} (cf.~\cite[
Theorem 4.8(b) and (c)]{GW89} or \cite[Theorem 3.30]{Rot79}): 
\bigskip

\noindent Let $M$ be an $R$-module and let $(I_R(M),\iota_M)$ be an injective hull of $M$. 

\begin{itemize}
\item[(I)] If $I$ is an injective $R$-module and $\iota$ is an $R$-module monomorphism from 
           $M$ into $I$, then every $R$-module homomorphism $\eta$ from $I_R(M)$ into $I$ 
           with $\eta\circ\iota_M=\iota$ is a monomorphism. (Since $I$ is injective and 
           $\iota_M$ is an $R$-module monomorphism, there {\it always\/} exists an $R$-module 
           homomorphism from $I_R(M)$ into $I$ with $\eta\circ\iota_M=\iota$\,!)
\item[(E)] If $N$ is an $R$-module and $\varphi$ is an $R$-module monomorphism from $M$ into 
           $N$ such that $\varphi(M)$ is an essential submodule of $N$, then every $R$-module 
           homomorphism $\nu$ from $N$ into $I_R(M)$ with $\nu\circ\varphi=\iota_M$ is a 
           monomorphism. (Since $I_R(M)$ is injective and $\varphi$ is an $R$-module monomorphism, 
           there {\it always\/} exists an $R$-module homomorphism from $N$ into $I_R(M)$ with 
           $\nu\circ\varphi=\iota_M$\,!)
\end{itemize}

\noindent (I) says that injective hulls are {\it minimal injective extensions\/} and (E) says 
that injective hulls are {\it maximal essential extensions\/}. In particular, injective hulls 
are uniquely determined up to isomorphism (cf.~\cite[Proposition 4.9]{GW89}).

Recall that a module is said to be {\it locally finite\/} if every finitely generated (or
equivalently, every cyclic) submodule is finite-dimensional. 

\begin{thm}\label{ess}
Let $\ag$ be a finite-dimensional Lie algebra over a field of prime characteristic. Then every 
essential extension of a locally finite $\ag$-module is locally finite.
\end{thm}

\noindent {\it Proof.\/} Let $M$ be a locally finite $\ag$-module, let $E$ be an essential 
extension of $M$, and let $e$ be any non-zero element of $E$. Then $E^\prime:=U(\ag)e$ is 
an essential extension of $M^\prime:=E^\prime\cap M$. Since $U(\ag)$ is noetherian, $M^\prime
\subseteq E^\prime$ is finitely generated. Because $M$ is by assumption locally finite, 
$M^\prime\subseteq M$ is finite-dimensional. By virtue of Theorem \ref{FBN}, we can apply 
\cite[Corollary 3.6]{Jat74} or the main result of \cite{Sch75} which both show that $E^\prime$ 
is also finite-dimensional, i.e., $E$ is locally finite.\quad $\Box$
\bigskip

\noindent The next result is an immediate consequence of Theorem \ref{ess}. 

\begin{cor}\label{locfincharp}
If $\ag$ is a finite-dimensional Lie algebra over a field of prime characteristic, then the 
injective hull of every locally finite $\ag$-module is locally finite.\quad $\Box$
\end{cor}

It is well-known that Corollary \ref{locfincharp} is also true for a finite-dimensional 
{\it solvable\/} Lie algebra over an arbitrary field of characteristic {\it zero\/} (see 
\cite[Theorem 2.2.3]{Don82} and \cite[Corollary 12]{Dah84}), but it does {\it not\/} hold 
for a finite-dimensional {\it semisimple\/} Lie algebra over a field of characteristic 
{\it zero\/} (see \cite[Remark after the proof of Proposition 2.2.2]{Don82} and \cite[Remark 
1]{Dah89}). More precisely, we have the following result.

\begin{thm}\label{locfinchar0}
Let $\ag$ be a finite-dimensional Lie algebra over a field of characteristic zero. Then the 
following statements are equivalent:
\begin{enumerate}
\item[{\rm (1)}] $\ag$ is solvable.
\item[{\rm (2)}] The injective hull of the one-dimensional trivial $\ag$-module is locally 
                 finite.
\item[{\rm (3)}] The injective hull of every locally finite $\ag$-module is locally finite.
\end{enumerate}
\end{thm}

\noindent {\it Proof.\/} The implication (1)$\Longrightarrow$(3) is just \cite[Theorem 
2.2.3]{Don82} or \cite[Corollary 12]{Dah84} and the implication (3)$\Longrightarrow$(2) 
is trivial. Hence it only remains to show the implication (2)$\Longrightarrow$(1).

Suppose that the injective hull $I_\ag(\fb)$ of the one-dimensional trivial $\ag$-module 
$\fb$ is locally finite. Since the ground field is assumed to have characteristic zero, 
the Levi decomposition theorem (cf.~\cite[p.~91]{Jac79}) yields the existence of a semisimple 
subalgebra $\s$ of $\ag$ (a so-called {\it Levi factor\/} of $\ag$) such that $\ag$ is the 
semidirect product of $\s$ and its solvable radical $\sol(\ag)$. According to \cite[Proposition
4]{Dah84}, the restriction $I:=I_\ag(\fb)_{\mid\s}$ is an injective $U(\s)$-module. Since
$I_\ag(\fb)$ is a locally finite $U(\ag)$-module, $I$ is a locally finite $U(\s)$-module.

Since $I$ is injective, it follows from the universal property (I) of injective hulls that
$\fb\subseteq I_\s(\fb)\subseteq I$. If $0\ne m\in I_\s(\fb)$, then the cyclic submodule 
$M:=U(\s)m$ of $I$ is finite-dimensional. Since $I_\s(\fb)$ is an essential extension of
$\fb$ and $M$ is a non-zero submodule of $I_\s(\fb)$, $M\cap\fb\ne 0$. Then for dimension
reasons, $M\cap\fb=\fb$, i.e., $\fb\subseteq M$. By virtue of Weyl's completely 
reducibility theorem (cf.~\cite[Theorem III.8, p.~79]{Jac79}), $\fb$ has a complement 
in $M$, i.e., there exists a submodule $C$ of $M$ such that $M=\fb\oplus C$. In particular, 
$\fb\cap C=0$ which implies that $C=0$ because $C$ is a submodule of $I_\s(\fb)$. Consequently, 
$M=\fb$ and therefore $\fb=I_\s(\fb)$. Hence $\fb$ is an injective $U(\s)$-module and thus 
also an injective $U(\fb s)$-module for every element $s\in\s$ (cf.~\cite[Proposition 
4]{Dah84}). Finally $\ext_{U(\fb s)}^1(\fb,\fb)\cong H^1(\fb s,\fb)\ne 0$ for every $0\ne 
s\in\s$ yields $\s=0$, i.e., $\ag=\sol(\ag)$ is solvable.\quad $\Box$
\bigskip

Let $\ag$ be a finite-dimensional Lie algebra over a field of characteristic zero. Donkin
\cite[Theorem 2.2.3]{Don82} proved that the largest locally finite submodule $I_\ag(M)_{{
\rm loc}}$ of the injective hull of any {\it finite-dimensional\/} $\ag$-module $M$ is 
{\it artinian\/}. In particular, if $\ag$ is solvable, then injective hulls of {\it 
finite-dimensional\/} $\ag$-modules are {\it artinian\/}. Furthermore, Dahlberg \cite{Dah89} 
showed that the injective hull of every {\it artinian\/} $\slg_2(\cb)$-module is {\it 
locally artinian\/}. In prime characteristic the following stronger result holds.

\begin{thm}\label{art}
If $\ag$ is a finite-dimensional Lie algebra over a field of prime characteristic, then the 
injective hull of every artinian $\ag$-module is artinian.
\end{thm}

\noindent {\it Proof.\/} Let $M$ be an artinian $\ag$-module. Then the socle $\soc_\ag(M)$ of 
$M$ is also artinian, i.e., a {\it finite\/} direct sum of simple modules. According to $I_\ag(M)
\cong I_\ag(\soc_\ag(M))$ and the additivity of $I_\ag(-)$, the assertion is an immediate 
consequence of (PI) and \cite[Theorem 2]{Jat75}.\quad $\Box$
\bigskip

Non-zero {\it noetherian\/} $\ag$-modules are very often {\it not injective\/}. This was proved 
in \cite[Corollary 2.3]{BHM82} for every (not necessarily commutative) {\it local\/} noetherian 
associative ring and motivated the first part of Proposition \ref{injdim} below. In particular, 
injective hulls of noetherian (or even finite-dimensional) $\ag$-modules are not noetherian. 
Moreover, for artinian and locally finite $\ag$-modules any possible injective dimension can 
occur. 

\begin{pro}\label{injdim}
Let $\ag$ be a finite-dimensional Lie algebra over a field $\fb$ of prime characteristic. Then
the following statements hold:
\begin{enumerate}
\item[{\rm (1)}] For every non-zero finitely generated (= noetherian) $\ag$-module $M$, 
                 we have $$\id_{U(\ag)}M=\dim_\fb\ag.$$
\item[{\rm (2)}] For every integer $0\le r\le\dim_\fb\ag$ there exists an artinian 
                 $\ag$-module $M_r$ such that $$\id_{U(\ag)}M_r=r.$$
\item[{\rm (3)}] For every integer $0\le r\le\dim_\fb\ag$ there exists a locally finite
                 $\ag$-module $N_r$ such that $$\id_{U(\ag)}N_r=r.$$
\end{enumerate}
\end{pro}

\noindent {\it Proof.\/} (1): Since $M$ is noetherian, it has a maximal submodule $N$.
Hence $S:=M/N$ is simple and thus finite-dimensional (cf.~\cite[Theorem 5.2.4]{SF88}). 
By virtue of \cite[Theorem 4.2(3)]{Far90}, there exists an $\ag$-module $V$ such that 
$\ext_{U(\ag)}^d(V,S)\ne 0$, where $d:=\dim_\fb\ag$. Then the long exact cohomology 
sequence implies the exactness of $$\ext_{U(\ag)}^d(V,M)\longrightarrow\ext_{U(\ag)}^d
(V,S)\longrightarrow\ext_{U(\ag)}^{d+1}(V,N).$$ Because of $\gl~U(\ag)=d$ (cf.~\cite[
Theorem 8.2]{CE56}), the right-hand term vanishes. One concludes that $\ext_{U(\ag)}^d
(V,M)\ne 0$, i.e., $\id_{U(\ag)}M\ge d$. The reverse inequality follows from $\id_{U(\ag)}
M\le\gl~U(\ag)=d$.

(2): Put $d:=\dim_\fb\ag$ and let $M_d$ be any non-zero finite-dimensional $\ag$-module. 
By the first part, we have $\id_{U(\ag)}M_d=d$. According to Theorem \ref{art}, the 
injective hull $I_\ag(M_d)$ and therefore $M_{d-1}:=I_\ag(M_d)/M_d$ are artinian. 
From the long exact cohomology sequence and the injectivity of $I_\ag(M_d)$ one 
concludes for an arbitrary $\ag$-module $X$ that $$\ext_{U(\ag)}^d(X,M_{d-1})\cong
\ext_{U(\ag)}^{d+1}(X,M_d)=0$$ because $\id_{U(\ag)}M_d=d$. Hence $\id_{U(\ag)}
M_{d-1}\le d-1$ (cf.~\cite[Theorem 9.8]{Rot79}). By another application of 
\cite[Theorem 9.8]{Rot79}, there exists an $\ag$-module $X_d$ such that $\ext_{U(\ag)}^d
(X_d,M_d)\ne 0$. Then the long exact cohomology sequence implies $$\ext_{U(\ag)}^{d-1}
(X_d,M_{d-1})\cong\ext_{U(\ag)}^d(X_d,M_d)\ne 0,$$ i.e., $\id_{U(\ag)}M_{d-1}=d-1$, 
and the assertion follows by induction.

(3): The proof is the same as for (2) except that one uses Corollary \ref{locfincharp}
instead of Theorem \ref{art} to conclude that $N_{d-1}:=I_\ag(N_d)/N_d$ is locally finite. 
\quad $\Box$
\bigskip

\noindent {\bf Remark.} Dually, non-zero {\it artinian\/} $\ag$-modules are {\it never 
projective\/} if $\ag\neq 0$ and for noetherian $\ag$-modules any possible projective 
dimension can occur (see \cite{Fel90}).
\bigskip

\noindent Since every {\it simple\/} module is {\it finitely generated\/}, the following 
is an immediate consequence of Proposition \ref{injdim}(1). 

\begin{cor}\label{injdimsim}
Let $\ag$ be a finite-dimensional Lie algebra over a field $\fb$ of prime characteristic 
and let $S$ be a simple $\ag$-module. Then $\id_{U(\ag)}S=\dim_\fb\ag.$\quad $\Box$
\end{cor}


\section{Locally Finite Submodules of the Coregular Module}


Let $\ag$ be a Lie algebra over a field $\fb$ of arbitrary characteristic. Then
the {\it linear dual\/} $U(\ag)^*:=\Hom_\fb(U(\ag),\fb)$ of $U(\ag)$ is a left
and a right $U(\ag)$-module, the so-called {\it coregular module\/} of $U(\ag)$ 
(cf.~\cite[2.7.7]{Dix96}). It is well-known that $U(\ag)^*$ is injective as a 
left and right $U(\ag)$-module (cf.~\cite[Proposition 1]{Lev76}).

Let $U(\ag)^\circ$ denote the {\it continuous dual\/} of $U(\ag)$ which is the 
largest locally finite submodule of the left and right $U(\ag)$-module $U(\ag)^*$. 
It is well-known that $U(\ag)^\circ$ also consists of all linear forms on $U(\ag)$ 
that vanish on some two-sided ideal of finite codimension in $U(\ag)$ 
(cf.~\cite[p.~51]{Lev86}). 

Finally, let $U(\ag)^\natural$ denote the set of all linear forms on $U(\ag)$ that
vanish on a certain power of the augmentation ideal $U(\ag)^+$ of $U(\ag)$. Then
one has the following inclusions where $\fb^*$ is identified with the linear forms
on $U(\ag)$ that vanish on $U(\ag)^+$ (cf.~\cite[Lemma 2.5.1]{Dix96}): $$\fb\cong\fb^*
\subseteq U(\ag)^\natural\subseteq U(\ag)^\circ\subseteq U(\ag)^*\,.$$ 

The following is also well-known (cf.~\cite[Lemme 2]{Lev76}).

\begin{lem}\label{essdual}
If $\ag$ is a finite-dimensional Lie algebra over an arbitrary field, then 
$U(\ag)^\natural$ is an essential extension of the one-dimensional trivial 
left and right $U(\ag)$-module.\quad $\Box$
\end{lem}

For the convenience of the reader we include a proof of the following result.

\begin{thm}\label{injnilpdual}\mbox{\rm (cf.~\cite[Th\'eor\`eme 3]{Lev76} or \cite[
Theorem 3]{Dah84})} If $\ag$ is a finite-dimensional nilpotent Lie algebra over an 
arbitrary field, then $U(\ag)^\natural$ is an injective hull of the one-dimensional 
trivial left and right $U(\ag)$-module.
\end{thm}

\noindent {\it Proof.\/} Since $\fb\cong\fb^*\subseteq U(\ag)^*$ and $U(\ag)^*$ is
injective, the universal property (I) of injective hulls implies that $I_\ag(\fb)
\subseteq U(\ag)^*$. It follows from \cite[Proposition 1]{Dah84} (cf.~also Theorem 
\ref{locfinchar0} and Corollary \ref{locfincharp}) that $I_\ag(\fb)$ is locally 
finite. Consider $\varphi\in I_\ag(\fb)$. Then $E:=U(\ag)\varphi$ is a 
finite-dimensional extension of $\fb$. An application of Fitting's lemma 
(cf.~\cite[Theorem II.4, p.~39]{Jac79}) shows that $\ag$ acts nilpotently on $E$ and 
it follows from the Engel-Jacobson theorem (cf.~\cite[Corollary 1.3.2]{SF88}) that a 
certain power of the augmentation ideal $U(\ag)^+$ annihilates $E$. Consequently, 
$\varphi\in U(\ag)^\natural$ and therefore $I_\ag(\fb)\subseteq U(\ag)^\natural$. 
Finally, the other inclusion follows from Lemma \ref{essdual} and the universal 
property (E) of injective hulls.\quad $\Box$
\bigskip

\noindent {\bf Remark.} It is observed in \cite[Remarque 2 after Th\'eor\`eme 3]{Lev76}
that $U(\ag)^\natural$ is {\it not\/} injective for the two-dimensional non-nilpotent
Lie algebra. It would be interesting to know whether the injectivity of $U(\ag)^\natural$ 
implies that $\ag$ is nilpotent.
\bigskip

\noindent The isomorphism $H^n(\ag,U(\ag)^\natural)\cong\ext_{U(\ag)}^n(\fb,U(\ag)^\natural)$ 
in conjunction with Theorem \ref{injnilpdual} and \cite[Theorem 7.6]{Rot79} yields the following 
cohomological vanishing theorem due to Koszul:

\begin{cor}\label{nilpcohvan}\mbox{\rm (cf.~\cite[Th\'eor\`eme 6]{Kos54})}
If $\ag$ is a finite-dimensional nilpotent Lie algebra over an arbitrary field, then 
$$H^n(\ag,U(\ag)^\natural)=0$$ for every positive integer $n$.\quad $\Box$
\end{cor}

\noindent {\bf Question.} Does the vanishing $H^n(\ag,U(\ag)^\natural)$ for every positive 
integer $n$ imply that $\ag$ is nilpotent? 
\bigskip

Let us now consider arbitrary finite-dimensional Lie algebras over fields of prime characteristic.

\begin{thm}\label{injcontdualcharp}
If $\ag$ is a finite-dimensional Lie algebra over a field of prime characteristic, then 
the continuous dual $U(\ag)^\circ$ is injective as a left and right $U(\ag)$-module.
\end{thm}

\noindent {\it Proof.\/} Since $U(\ag)^\circ\subseteq U(\ag)^*$ and $U(\ag)^*$ is
injective, the universal property (I) of injective hulls implies that $I_\ag(U(\ag)^\circ)
\subseteq U(\ag)^*$. Because $U(\ag)^\circ$ is locally finite, it follows from Corollary 
\ref{locfincharp} that $I_\ag(U(\ag)^\circ)$ is also locally finite. But since by 
definition $U(\ag)^\circ$ is the {\it largest\/} locally finite submodule of $U(\ag)^*$,
$U(\ag)^\circ=I_\ag(U(\ag)^\circ)$ is injective.\quad $\Box$
\bigskip

\noindent The isomorphism $H^n(\ag,U(\ag)^\circ)\cong\ext_{U(\ag)}^n(\fb,U(\ag)^\circ)$ 
in conjunction with Theorem \ref{injcontdualcharp} and \cite[Theorem 7.6]{Rot79} yields 
the following cohomological vanishing theorem:

\begin{cor}\label{cohvan}
If $\ag$ is a finite-dimensional Lie algebra over a field of prime characteristic, then 
$$H^n(\ag,U(\ag)^\circ)=0$$ for every positive integer $n$.\quad $\Box$
\end{cor}

\noindent {\bf Remark.} The case $n=1$ of Corollary \ref{cohvan} was already proved by Masuoka 
\cite[Proposition 5.1]{Mas00}. It follows from Corollary \ref{cohvan} in conjunction with
\cite[Th\'eor\`eme 2]{Kos54} that every cohomology class of a finite-dimensional Lie algebra 
over a field of prime characteristic with coefficients in a finite-dimensional module is
annihilable. This result was proved in a completely different way by Dzhumadil'daev \cite[Theorem 
3.1, pp.~467--470]{Dzh90}.
\bigskip

The equivalence of (1), (3), and (4) in the next result is essentially due to Koszul (see 
\cite[Th\'eor\`eme 7 and p.~536]{Kos54}. Moreover, for an {\it algebraically closed\/} 
ground field the implication (1)$\Longrightarrow$(2) follows from \cite[Th\'eor\`eme 3.6 
and Th\'eor\`eme 4.10]{BM85} (see also \cite[Proposition 3.4 and Proposition 3.6]{Lev86} 
for $\fb=\cb$).

\begin{thm}\label{injcontdualchar0}
Let $\ag$ be a finite-dimensional Lie algebra over a field of characteristic zero. Then the 
following statements are equivalent:
\begin{enumerate}
\item[{\rm (1)}] $\ag$ is solvable.
\item[{\rm (2)}] The continuous dual $U(\ag)^\circ$ is injective as a left and right 
                 $U(\ag)$-module.
\item[{\rm (3)}] $H^n(\ag,U(\ag)^\circ)=0$ for every positive integer $n$.
\item[{\rm (4)}] $H^3(\ag,U(\ag)^\circ)=0$.
\end{enumerate}
\end{thm}

\noindent {\it Proof.\/} The proof of the implication (1)$\Longrightarrow$(2) is the 
same as for Theorem \ref{injcontdualcharp} except that one uses Theorem \ref{locfinchar0}
instead of Corollary \ref{locfincharp} in order to conclude that $I_\ag(U(\ag)^\circ)$ 
is locally finite. Since (2)$\Longrightarrow$(3) is clear and (4) is just a special case 
of (3), it remains to show the implication (4)$\Longrightarrow$(1).

Suppose that $H^3(\ag,U(\ag)^\circ)=0$ and let $M$ be an arbitrary finite-dimensional 
$\ag$-module. Then the isomorphism $\Hom_\fb(U(\ag),M)\cong U(\ag)^*\otimes_\fb M$ 
(where $M$ is considered as a {\it trivial\/} $\ag$-module) implies that $H^3(\ag,
\Hom_\fb(U(\ag),M)_\loc)=0$ where $\Hom_\fb(U(\ag),M)_\loc$ denotes the largest 
locally finite submodule of $\Hom_\fb(U(\ag),M)$. According to \cite[Th\'eor\`eme 2]{Kos54}, 
it follows that every cohomology class in $H^3(\ag,M)$ is annihilable and thus \cite[5), 
p.~536]{Kos54} yields that $\ag$ is solvable.\quad $\Box$
\bigskip

\noindent {\bf Remark.} The above proof of the implication (1)$\Longrightarrow$(2) is not 
only much more direct than in \cite{BM85} or \cite{Lev86} but also answers affirmatively 
a question posed at the end of the third section in \cite{BM85}. Moreover, it should be 
noted that the implication (2)$\Longrightarrow$(1) in Theorem \ref{injcontdualchar0} can 
also be obtained directly from the universal property (I) of injective hulls and Theorem 
\ref{locfinchar0}. 
\bigskip

\noindent Recently, H.-J. Schneider has generalized the implication (1)$\Longrightarrow$(2) 
in Theorem \ref{injcontdualchar0} even further. Let $\ag$ be a finite-dimensional Lie
algebra over a field of characteristic zero and let $\sol(\ag)$ denote the solvable 
radical of $\ag$. Then Schneider proves that the restriction $[U(\ag)^\circ]_{\mid
\sol(\ag)}$ of $U(\ag)^\circ$ to $\sol(\ag)$ is injective (cf.~\cite[Theorem 5.3]{Mas00}). 
This in conjunction with the Hochschild-Serre spectral sequence (cf.~\cite[Theorem 6]{HS53}) 
and the two Whitehead lemmata (cf.~\cite[Theorem III.13]{Jac79}) implies that 
$H^1(\ag,U(\ag)^\circ)=0=H^2(\ag,U(\ag)^\circ)$ (see \cite[Proposition 5.1 and Theorem 
5.2]{Mas00}). But Theorem \ref{injcontdualchar0} shows that $H^3(\ag,U(\ag)^\circ)\ne 0$ 
if $\ag$ is {\it not solvable\/} which generalizes \cite[Remark 5.9]{Mas00}.
\bigskip

It follows from the universal properties (E) and (I) of injective hulls in conjunction with
Lemma \ref{essdual} and Theorem \ref{injcontdualcharp} or Theorem \ref{injcontdualchar0} that 
$$\fb\cong\fb^*\subseteq U(\ag)^\natural\subseteq I_\ag(\fb)\subseteq U(\ag)^\circ\,.$$

\noindent Note that the cocommutative Hopf algebra structure on $U(\ag)$ induces a commutative
algebra structure on $U(\ag)^*$ which over a field $\fb$ of characteristic zero can be identified 
with the algebra of power series in $\dim_\fb\ag$ variables (cf.~\cite[Proposition 2.7.5]{Dix96}) 
and the continuous dual $U(\ag)^\circ$ is a subalgebra of $U(\ag)^*$. 

Let $\ag$ be a finite-dimensional {\it solvable\/} Lie algebra over the complex numbers. Then 
Levasseur \cite[Th\'eor\`eme 2.2]{Lev86} has shown that $I_\ag(\fb)$ is isomorphic to a polynomial 
algebra in $\dim_\fb\ag$ variables on which $\ag$ acts via derivations.
\bigskip

\noindent {\bf Conjecture.} Let $\ag$ be a finite-dimensional Lie algebra over a field $\fb$. If 
$\ag$ is solvable and char$(\fb)=0$ or if $\ag$ is arbitrary and char$(\fb)>0$, then $I_\ag(\fb)$ 
is isomorphic to a polynomial algebra in $\dim_\fb\ag$ variables on which $\ag$ acts via derivations.
\bigskip

\noindent If $\ag$ is abelian, then this follows from \cite[Theorem 2]{Nor74} and in 
\cite[Section 4]{Dah84} there are examples confirming this for Lie algebras of small 
dimensions.


\section{Minimal Injective Resolutions}


Let $\ic$ be a two-sided ideal of an associative ring $R$. Then $$\ug(\ic):=\sup\{n\in\nb_0
\mid\ext_R^n(R/\ic,R)\ne 0\}$$ and $$\lgr(\ic):=\inf\{n\in\nb_0\mid \ext_R^n(R/\ic,R)\ne 0\}$$ 
are called {\it upper grade\/} and {\it lower\/} (or {\it homological\/}) {\it grade\/} of $\ic$, 
respectively. A left and right noetherian associative ring $R$ is left (resp.~right) {\it 
injectively homogeneous\/} over a central subring $C$ if $R$ is integral over $C$, $\id_R 
R<\infty$ (resp.~$\id R_R<\infty$) and $\ug(\mc)=\ug(\nc)$ for all maximal ideals $\mc$ and 
$\nc$ such that $\mc\cap C=\nc\cap C$. In \cite{BH88} it was demonstrated that for associative
rings integral over a central subring the class of injectively homogeneous rings is a natural 
generalization of the class of commutative {\it Gorenstein rings\/}. Moreover, \cite[Corollary 
3.6]{BH88} shows that $R$ is injectively homogeneous over its center $C(R)$ if and only if 
$R$ is injectively homogeneous over {\it every\/} subring $C\subseteq C(R)$ over which $R$ 
is integral, and by virtue of \cite[Corollary 4.4]{BH88}, $R$ is left injectively homogeneous 
if and only if $R$ is right injectively homogeneous.

\begin{lem}\label{injhom}
If $\ag$ is a finite-dimensional Lie algebra over a field of prime characteristic, then 
$U(\ag)$ is injectively homogeneous over every subring $C$ of $U(\ag)$ with $\oc(\ag)
\subseteq C\subseteq C(U(\ag))$.
\end{lem}

\noindent {\it Proof.\/} Let $\mc$ be a maximal ideal of $U(\ag)$. Then $\wp:=\mc\cap 
C(U(\ag))$ is also maximal \cite[Corollary 6.3.4]{SF88}, and thus Hilbert's Nullstellensatz 
yields that $C(U(\ag))/\wp$ is finite-dimensional. Since $U(\ag)$ is finitely generated 
over $C(U(\ag))$, we conclude that $M:=U(\ag)/\mc$ is also finite-dimensional. Set $d:=
\dim_\fb\ag$. According to \cite[Theorem 4.2(3)]{Far90}, there exists a simple $\ag$-module 
$S$ such that $\ext_{U(\ag)}^d(M,S)\ne 0$. If $\ac$ denotes the annihilator of a generator 
of $S$ in $U(\ag)$, we obtain a short exact sequence $0\to\ac\to U(\ag)\to S\to 0$ of $U(\ag)
$-modules. The long exact cohomology sequence implies the exactness of $$\ext_{U(\ag)}^d(M,
U(\ag))\longrightarrow\ext_{U(\ag)}^d(M,S)\longrightarrow\ext_{U(\ag)}^{d+1}(M,\ac).$$ 
Because of $\gl~U(\ag)=d$ (cf.~\cite[Theorem 8.2]{CE56}), the right-hand term vanishes. We 
conclude that $\ext_{U(\ag)}^d(M,U(\ag))\ne 0$, i.e., $\ug(\mc)\ge\lgr(\mc)\ge d$. The 
reverse inequality follows from $\ug(\mc)\le\gl~U(\ag)=d$. Hence $\ug(\mc)=d$ for every 
maximal ideal of $U(\ag)$. This and (IC) yield the assertion.\quad $\Box$
\bigskip

\noindent {\bf Remark.} Let $\ag$ be a finite-dimensional Lie algebra over a field of  
characteristic zero. According to a theorem of Laty\v sev \cite{Lat63}, $U(\ag)$ is a PI 
algebra if and only if $\ag$ is abelian. Since every algebra which is a finitely generated 
module over its center is a PI algebra (cf.~\cite[p.~xi]{GW89}), $U(\ag)$ is injectively 
homogeneous over its center if and only if $\ag$ is abelian. 
\bigskip

One consequence of Lemma \ref{injhom} is that $\id\, U(\ag)_{\wp}<\infty$ for every 
semiprime ideal $\wp$ of every subring $C$ of $U(\ag)$ with $\oc(\ag)\subseteq C
\subseteq C(U(\ag))$ (cf.~\cite[Fundamental Theorem (e), p.~10]{Bas63} and 
\cite[Theorem 4.1]{BH88}). More importantly for the purpose of this paper, it is an 
immediate consequence of Lemma \ref{injhom} and \cite[Theorem 5.5]{BH88} that the minimal
injective resolution of $U(\ag)$ has the same form as for {\it commutative Gorenstein 
rings\/} (cf.~\cite[Fundamental Theorem (f), p.~10]{Bas63}). Recall that a {\it minimal 
injective resolution\/} of a module $M$ is a long exact sequence $$0\longrightarrow M
\longrightarrow I_0\stackrel{d_0}{\longrightarrow}I_1\longrightarrow\cdots\longrightarrow 
I_n\stackrel{d_n}{\longrightarrow}I_{n+1}\longrightarrow\cdots$$ such that $I_n$ is an 
injective hull of $\Ker(d_n)$ for every non-negative integer $n$.

\begin{thm}\label{mininjres}
Let $\ag$ be a finite-dimensional Lie algebra over a field $\fb$ of prime characteristic. 
If $0\longrightarrow U(\ag)\longrightarrow I_0\longrightarrow\cdots\longrightarrow I_d
\longrightarrow 0$ is a minimal injective resolution of $U(\ag)$ as a left or right 
$U(\ag)$-module, then $$I_n\cong\bigoplus_{\hoe(\pc)=n}I_\ag(U(\ag)/\pc)$$ for every 
$0\le n\le d:=\dim_\fb\ag$.\quad $\Box$
\end{thm}

\noindent {\bf Remark.} If $\ag$ is a finite-dimensional Lie algebra over a field of  
characteristic zero, then the structure of a minimal injective resolution of $U(\ag)$ 
is even in the {\it solvable\/} case more complicated than in Theorem \ref{mininjres} 
(cf.~\cite{Mal83,Mal86}). 
\bigskip

Let $\ag$ be a finite-dimensional solvable Lie algebra over an algebraically closed 
field of characteristic zero. Then the last term of a minimal injective resolution 
of $U(\ag)$ is isomorphic to the continuous dual $U(\ag)^\circ$ of $U(\ag)$ 
(see \cite[Th\'eor\`eme 3.6 and Th\'eor\`eme 4.10]{BM85} and also \cite[Proposition 
3.4 and Proposition 3.6]{Lev86} for $\fb=\cb$). We conclude the paper by applying 
Theorem \ref{mininjres} in order to prove the analogue of this result in prime 
characteristic.

\begin{thm}\label{lastermininjres}
Let $\ag$ be a finite-dimensional Lie algebra over an algebraically closed field of 
prime characteristic. If $0\longrightarrow U(\ag)\longrightarrow I_0\longrightarrow\cdots
\longrightarrow I_d\longrightarrow 0$ is a minimal injective resolution of $U(\ag)$ as 
a left or right $U(\ag)$-module, then $I_d\cong U(\ag)^\circ$.
\end{thm}

\noindent {\it Proof.\/} By virtue of Corollary \ref{locfincharp}, injective hulls of 
locally finite modules are locally finite. Since $\fb$ is algebraically closed, this 
enables one to prove that $$U(\ag)^\circ\cong\bigoplus_{S\in\irr(\ag)}I_\ag(S)^{\oplus
\dim_\fb S}$$ as a left or right $U(\ag)$-module, where $\irr(\ag)$ denotes the set of 
isomorphism classes of simple $\ag$-modules (cf.~\cite[1.5]{Gre76} for the analogous
statement in terms of coalgebras and comodules). On the other hand, it follows from 
(PI) and \cite[Theorem 4]{Sch76} that a prime ideal $\pc$ of $U(\ag)$ has maximal 
height $d$ if and only if $\pc$ is maximal. But every maximal ideal $\pc$ of $U(\ag)$ 
is primitive, i.e., there is a simple $\ag$-module $S$ such that $\pc=\ann_{U(\ag)}(S)$. 
Then the density theorem (cf.~\cite[Theorem 16, p.~95]{Kap72}) yields that 
$$U(\ag)/\pc=U(\ag)/\ann_{U(\ag)}(S)\cong\End_\fb(S)\cong S^{\oplus\dim_\fb S^*}$$ 
as a left or right $U(\ag)$-module. In particular, simple $\ag$-modules are isomorphic
if and only if their annihilators in $U(\ag)$ coincide. According to (PI) and Kaplansky's 
theorem (cf.~\cite[Theorem 50, p.~157]{Kap72}), every primitive ideal of $U(\ag)$ is maximal 
and therefore $$I_d\cong\bigoplus_{\hoe(\pc)=d}I_\ag(U(\ag)/\pc)\cong\bigoplus_{S\in
\irr(\ag)}I_\ag(S)^{\oplus\dim_\fb S^*}\cong U(\ag)^\circ\,.\quad\Box$$
\bigskip

\noindent {\bf Question.} Does Theorem \ref{lastermininjres} remain valid for {\it 
arbitrary\/} ground fields of prime characteristic?



\end{document}